# (1, k)-Swap Local Search for Maximum Clique Problem

LAVNIKEVICH N.

Given simple undirected graph G = (V, E), the Maximum Clique Problem(MCP) is that of finding a maximum-cardinality subset Q of V such that any two vertices in Q are adjacent. We present a modified local search algorithm for this problem. Our algorithm build some maximal solution and can determine in polynomial time if a maximal solution can be improved by replacing a single vertex with $k$, $k \geq 2$, others.

We test our algorithms on DIMACS[5], Sloane[15], BHOSLIB[1], Iovanella[8] and our random instances.

## *1. Introduction*

Given an arbitrary *simple(*without loops and multiple edges) *undirected graph G = (V, E)*. The *order* of a graph is number of vertices $n = |V|$, where $V = \{1, \ldots, n\}$ is *vertex set*.

Two vertices *v* and *u* of a graph *G* may be adjacent or not. If two vertices *v* and *u* are adjacent they perform an edge e*(v, u)*. All edges perform *edges set E*. A graph's *size* is the number of edges in the graph, $|E| = m$.

The neighbourhood *N(v)* of a vertex *v* in a graph G consists of all vertices adjacent to *v*, $N(v) = \{ u \mid u \in V, e(v, u) \in E \}$. By $\overline{N}(v)$ we will denote set of all vertices not adjacent with *v*, $\overline{N}(v) = \{ u \mid u \in V, u \neq v \text{ and } e(v, u) \notin E \}$.

A *complement* of graph *G* is the graph $\overline{G} = (V, \overline{E})$ where $\overline{E} = \{ (v, u) \mid v, u \in V, v \neq u \text{ and } (v, u) \notin E \}$.

Let $U \subset V$ be any subset of vertices of G. Then the *induced subgraph* G[U] is the graph whose vertex set is U and whose edge set consists of all of edges in E that have both endpoints in U.

Important types of induced subgraphs are *cliques* and *independent sets*.

A *clique Q* is a set of pairwise adjacent vertices of the graph. The *maximum clique problem* (MCP) is to find a clique of maximum cardinality in a graph *G*.

An *independent set (stable set, vertex packing) I* is a set of pairwise nonadjacent vertices of the graph. The *maximum independent set* (MIS) problem is to find an independent set of maximum cardinality in a graph *G*.

In particular, a clique (independent set) which is not properly contained in any other clique is called *maximal*. A maximal clique (independent set) with the maximum size is called a *maximum* clique (independent set).

It is easy to see that *Q* is a clique in a graph *G = (V, E)* if and only if *Q* is an independent set of $\overline{G}$. Thus, the *Maximum Clique Problem(MCP)* and the *Maximum Independent Set(MIS) Problem* are equivalent.

In addition, both problems are $\mathcal{NP}-$complete, Karp[9], which means that unless $\mathcal{P}= \mathcal{NP}$ there exists no algorithm that can solve this problems in time polynomial to the order of the input graph.

Since all known exact algorithms [12, 4, 13, 17] for these problems take exponential time, making large graphs infeasible to solve in practice. In [11] we showed that usage of coloring heuristics in exact MCP algorithms provides $\Omega(2^{0.2n})$ algorithm running time lower bound.

Instead, heuristic algorithms [1, 4, 6, 12, 16] are used to efficiently compute high-quality cliques(independent sets). They widely used in exact algorithms to obtain an initial solution [2, 14, 17].

In this paper we will consider local search algorithm only for one of two problems – *Maximum Clique Problem (MCP)*.

## 2. Heuristic Algorithms

There are a wide range of heuristics and local search algorithms for the MCP problem (see for example Pardalos [12], Bomse [4]). The majority of this approximation algorithms in the literature are called *sequential greedy heuristics*. These heuristics generate a maximal clique through the repeated *addition* of a vertex into a partial clique or the repeated *deletion* of a vertex from a set that is not a clique.

Kopf and Ruhe [10] named these two classes of heuristics *the Best in* and *the Worst out* heuristics. Decisions on which vertex to be added in or moved out next are based on certain indicators associated with candidate vertices. For example, a possible Best in heuristic constructs a maximal clique by repeatedly adding in a vertex that has the largest degree among candidate vertices. In this case, the indicator is the degree of a vertex. On the other hand, a possible Worst out heuristic can start with the whole vertex set V. It will repeatedly remove a vertex out of V until V becomes a clique.

Kopf and Ruhe [10] further divided the above two classes of heuristics into *New* and *Old* (Best in or Worst out) heuristics. Namely, if the indicators are updated every time a vertex is added in or moved out, then the heuristic is called a *New* heuristic. Otherwise it is called an *Old* heuristic. We can find in the literature that most heuristics for the maximum clique problem fall in one or the other classes. See for example, the approximation algorithms E1 and E2 of Johnson [6], and the approximation algorithm NMCLIQ of Tomita et al [16]. The differences among these heuristics are their choice of indicators and how indicators are updated. A heuristic of this type can run very fast.

A common feature of the *sequential heuristics* is that they all find only one maximal clique. Once a maximal clique is found, the search stops.

We can view this type of heuristics from a different point of view. Let us define $S_G$ to be the space consisting of all the maximal cliques of G. What a sequential greedy heuristic does is to find one point in $S_G$, hoping it is (close to) the optimal point. This suggests us a possible way to improve our approximation solutions, namely, expand the search in $S_G$. For example, once we find a point $x \in S_G$, we can search its neighbors to improve x. This leads to the class of the *local search heuristics*.

One of the best local search heuristics is *iterated local search(ILS)* algorithm by Andrade et al. [1]. This local search algorithm uses (1, 2)-*swaps* or 2-*improvements* to gradually increase the size of the current solution Q. This is done by taking one vertex *u* from the solution Q and inserting two another vertices into the solution. In general, a (j, k)-swap takes out j vertices and inserts k vertices.

A vertex is called *k-tight*, if exactly k of its not neighbors lie in the solution. A 2-improvement can be applied if the inserted adjacent vertices x, y are 1-tight vertices with common vertex v, which is to be removed from the solution Q.

A 0-tight vertex is called *free* and can be added to the solution without any restrictions. ILS stores the tightness for each vertex and updates it after each swap. The local search iterates over all vertices $v \in V$ and computes the set of 1-tight vertices which, by definition, have only one not adjacent vertex *u* in the solution Q. If there are at least two adjacent candidates *v1* and *v2* both not adjacent with only one vertex u ∈ Q then the 2-improvement can be applied. Afterwards, the tightness of *v*'s not neighbors is decreased by one and the tightness of the inserted vertices not neighbors are increased by one.

We slightly improve this procedure.

## 3. Andrade 2–improvement Local Search Algorithm

Proposed by Andrade[1] procedure can be described as follows (Algorithm 1) .

Step A1.1. Apply some sequential heuristic algorithm and generate some maximal initial solution *Q*.

Step A1.2. Build candidates sets and set T1 of 1-tight vertices.

Step A1.3. If in T1 exists two adjacent vertices *v1* and *v2* which both not adjacent exactly with the same vertex *u* ∈ Q goto Step A1.4.

Else stop.

Step A1.4. Delete vertex *u* from solution *Q*, insert vertices *v1* and *v2* into solution *Q*, update set of candidates and T1 set. Go to Step A1.3.

---

*Algorithm 1*: ILS(SG, Q)
Global : graph G;
Input:   SG  - subgraph(as set of vertices);
Output:  Q – maximal clique .
**-------------------------------------------------------------------------------**
1. Heuristic(SG, Q);  // initiate solution Q - apply one of sequential  heuristic algorithms;
// build candidates set C;
2. Build_candidate_sets_and_T1(SG, Q, C);
// iteration local search step
3. while (Exists_improvement(SG, Q, C, $u_{swap}$, $C_{swap}$))
4.    Improve_solution_and_modify_a_set_of_candidates (SG, Q, C, $u_{swap}$, $C_{swap}$);

---

## 4. Local Search (1, k)–swap Algorithm

In our case we use bits-set(similar to [14]) implementation of the adjacent matrix and other sets. So if we need to add the vertex to the candidates set or remove it from the candidates set we only perform set_bit_to_1(set_bit_to_0) procedures. And if we need to remove, for example, adjacent vertices from the bits-set we perform 64-bit parallel AND_NOT operation.

We modified Algorithm 1 steps.

---
**Algorithm 2**: *Build_candidate_sets_and_T1* ( *SG, Q, C*)
Global : graph *G*;
Input:  *SG* - subgraph(as set of vertices);
        *Q* – current solution(maximal clique);
Output: *C* – candidate sets matrix.
--------------------------------------------------------------------------
1. *for each* vertex $v \in SG$ initiate empty candidates set $C[v] \in C$;
2. *for each* vertex $v \in SG \setminus Q$ {
   2.1. $C[v] = \overline{N}(v) \cap Q$;
   2.2. *if* $|C[v]| == 1$ *then* add $v$ into $C[u]$, where $u$ is the unique member of $C[v]$;
}
---

On the Step A1.2 we call procedure *Build_candidate_sets_and_T1* (**Algorithm 2**). For each vertex $v$ outside solution($v \notin Q$) we build candidates set $C[v]$ of not adjacent solution members: $C[v] = \overline{N}(v) \cap Q$. If $|C[v]| == 1$ (i.e. vertex $v$ is *1-tight*), we add $v$ into $C[u]$.

On the step A1.3 we analyze each solution Q members candidate sets $C[u]$, $u \in$

---
**Algorithm 3**: *Exists_improvement* (*SG, Q, C, $u_{swapped}$, $C_{improve}$*)
Global : graph *G*;
Input:  *SG* - subgraph(as set of vertices);
        *Q* –current solution(maximal clique);
        *C* – candidates set.
Output: $u_{swapped}$ – vertex to be removed from the solution *Q*;
        $C_{improve}$ – vertices which can be added to the solution *Q* instead of $u_{swap}$, $|C_{improve}| > 1$.
--------------------------------------------------------------------------------
1. Clear $C_{improve}$ set;
2. *for each* vertex $u \in Q$ {
   2.1. If $|C[u]| > 1$ && $|C[u]| > |C_{improve}|$ {
        2.1.1. fvs_QE ($C[u], Q_{local}$); // find quasi-exact solution for subgraph $C[u]$
        2.1.2. If $|Q_{local}| > 1$ && $|Q_{local}| > C_{improve}$ {
               2.1.2.1. $u_{swapped} = u$;
               2.1.2.2. Copy $Q_{local}$ to $C_{improve}$;
               2.1.2.3. }
        2.1.3. }
   2.2. }
3. *return* $|C_{improve}| > 1$
---

*Q*.

If some candidates set $C_u$ has more than one member it means that we can try to find improvement of the current solution *Q*.

Andrade[1] find 2-improvement. But we will try to build k-improvement, k ≥ 2.

---

**Algorithm 4**: *Improve_solution_and_modify_a_set_of_candidates* (SG, Q, C, $u_{swapped}$, $C_{improve}$))
*Global* : graph G;
*Input*: SG - subgraph(as set of vertices);
Q –current solution(maximal clique);
C – candidates set.
*Output*: $u_{swapped}$ – vertex to be removed from the solution Q;
$C_{improve}$ – vertices which can be added to the solution Q instead of $u_{swap}$, |$C_{improve}$| > 1.
-------------------------------------------------------------------------------
1. If |$C_{improve}$| < 2 return;
2. Remove $u_{swapped}$ from current solution Q;
3. Create tight candidates set $C[u_{improve}]$ : copy $C_{improve}$ to $C[u_{improve}]$;
4. For each $u \in C_{improve}$ clear candidates set $C_u$: add u into solution Q;
5. For each vertex $v \in SG \setminus Q \setminus u_{swapped}$ {
    5.1. remove $u_{swapped}$ and all not adjacent members of $C_{improve}$ from C[v];
    5.2. if |C[v]| == 1 then add v into C[u], where u is the unique member of C[v];
6. }

---

For this purpose we can use any heuristic maximal clique algorithm. But, taking into account the fact that usually the size of the candidates set C[u] is not large, we developed quasi-exact maximal clique algorithm(Algorithm fvs_QE) . We will describe it in next section.

On the Step A1.4, if |$C_{improve}$| < 2 we stop procedure and Q is maximal clique. Else(|$C_{improve}$| ≥ 2) we update solution and candidate sets: we remove $u_{swapped}$ from the current solution Q and create new candidates set $C[u_{swapped}]$. It is obvious that vertex $u_{swapped}$ is not adjacent only with $C_{improve}$ members and so new $Cu_{swapped} = C_{improve}$. Than for each vertex $v \in SG \setminus Q \setminus C_{improve} \setminus u_{swapped}$ we remove from candidates set C[v] vertex $u_{swapped}$ and all not adjacent members $C_{improve}$. If |C[v]| == 1 we add v into C[u] ,where u is the unique member of C[v]. Than we repeat step A1.3.

## 5. Quasi-exact Maximal Clique Algorithm

The main purpose of the quasi-exact maximal clique algorithm fvs_QE (Algorithm 5) is to find exact solution Q in a small N-vertex graphs. Under small graph we will understand graphs of order |G| ≤ N. For class of small graphs current algorithm will build exact solution in linear time. But we can apply this algorithm for

any undirected graph *G*. In this case algorithm will build not exact solution. So we

---

**Algorithm 5**: *Exists_ improvement* (SG, Q)
*Global* : graph *G*;
*Input*:   *SG*  - subgraph(as set of vertices);
*Output*: *Q* –current solution(maximal clique);
--------------------------------------------------------------------------------
1. Clear solution Q;
2. While (input subgraph SG is not empty) {
    2.1.1. Select from SG first  N(or less iff  |SG| < N, N – is constant currently == 6) vertices, remove each selected vertex from SG  and add it(vertex) into the current vertex set $SG_N$ ;
    2.1.2. Find exact solution in small graph  fvs_exact_in small_N_order_graph($SG_N$, $Q_N$);
    2.1.3. Q = Q ∪ $Q_N$;
    2.1.4. Remove from the input subgraph SG all vertices not adjacent with at least one vertex in the  local solution  $Q_N$.
3. }

---

called it as  "quasi-exact".

Exact maximum clique algorithm for small graphs currently use  N≤ 6. First of all, we check if $SG_N$ is empty. In this case we return empty solution. Otherwise for input subgraph $SG_N$ we build new adjacent matrix which we represent in  packed triangle form:

*(5,4), (5,3), (5,2), (5,1), (5,0), (4,3), (4,2), (4,1), (4,0), (3,2), (3,1), (3,0), (2,1), (2,0), (1,0).*

The length of this word must be *N*(N-1)/2* bits which we will number from 0 to *N*(N-1)/2-1*  from right to left. As we use N=6, we need 15-bit word to store adjacent matrix. Bit 0 corresponds to the graph edge (1,0), bit 1 – to the edge (2,0),  bit 2 – to the edge (2,1),  bit 10 – to the edge  (5,0),  and so on,  to the bit 15 which corresponds to the edge *(5,4)* .

Example 1.  Cycle $C_5$,  *V($C_5$)* = { 0, 1, 2, 3, 4 }, *E($C_5$)*={ (0,1), (1,2), (2,3), (3,4), (4,0) } has adjacent matrix  packed in   the  { 00000 1001 100 10 1 } 15-bits word.

It is known [7] that the number of labeled n-vertex simple undirected graphs is $2^{n*(n-1)/2}$.   So if we consider 6-vertex graph we need 32768 6-bit words. So we create CLIQUE_MEMBERS_MASK   array    of    32768    packed    solutions.    Each CLIQUE_MEMBERS_MASK   item is predefined packed clique items where in each packed solution byte we use bits from 0 to 5. If *i*-th bit set to  '1'  it means that *i*-th vertex of the input subgraph belongs to the clique.

Example 2. Let we have small 6-vertex subgraph $SG_N$ = $SG_6$ = { 14, 2, 138, 29, 77, 0 } and *E($SG_6$)*={ (0, 2), (0,29), (2,77), (2,138), (14,29) }  . We enumerate each $SG_6$ member by their position in subgraph, i.e. vertex  14 has 0-th index, vertex  2 has 1-th index, …, vertex 77 has 4-th index, vertex 0 has 5-th index. In accordance with new vertex indexes and graph G adjacent matrix we build packed 15-bits adjacent matrix $AM_{packed}$ = { 01010 0110 001 10 0 } and get result(packed clique members) = 010110 from CLIQUE_MEMBERS_MASK[$AM_{packed}$]. It means that items with indexes {1, 2, 4} perform the clique(remember that we index vertices of the subgraph $SG_6$ from left to right). In accordance with our enumeration of the input subgraph $SG_6$ we have that 010110 → { 1, 2, 4 } → { 2, 138, 77 }, i. e. vertices {2, 138, 77} perform maximum clique in $SG_N$ = { 14, 2, 138, 29, 77, 0 }.

Example 3. Let input subgraph $SG_N$ = $C_5$. We build packed triangle form $AM_{packed}$ = { 00000 1001 100 10 1 }. This bits vector corresponds to the 0x265 or 613 value. Predefined array CLIQUE_MEMBERS_MASK[613] contains packed clique members 000011 → {0, 1}.

Example 4. Let input subgraph $SG_N$ is complete graph $K_6$ = {187, 1, 2, 33, 5, 99 }. We build packed triangle form $AM_{packed}$ = { 11111 1111 111 11 1 }. This bits vector corresponds to the 0x7fff or 32767 value. Predefined array CLIQUE_MEMBERS_MASK[32767] contains packed clique members 111111 = {187, 1, 2, 33, 5, 99 }.

Example 5. Let input subgraph $SG_N$ is edgeless graph $I_2$ = { 65, 2 }. We build packed triangle form $AM_{packed}$ = { 00000 0000 000 00 0 }. This bits vector corresponds to the 0x0 or 0 value. Predefined array CLIQUE_MEMBERS_MASK[0] contains packed clique members 000001 → {0} → { 65 }.

## 6. Experimental Results

### 6.1. Measures

Our experiment (solution of the MCP or/and MIS problems) consists of
1. a measure defined for solved problem $\mathcal{P}$. In our case this is the size of maximal(maximum) clique $Q$.
2. a set of algorithms $\mathcal{A}_\mathcal{P}$;
3. a set of possible inputs (instances) $\mathcal{I}_\mathcal{P}$;
4. a map $\mathcal{S}_\mathcal{P}$ which for each $u \in \mathcal{I}_\mathcal{P}$ maps a finite set of algorithms $A \in \mathcal{A}_\mathcal{P}$ solutions.

As our goal is to find solution with maximal measure then we can specify our experiment as *maximization* problem.

The performance of a set of algorithms $\mathcal{A}$ for input $u \in \mathcal{I}_\mathcal{P}$ we define as

$$\mathcal{A}_\mathcal{P}(u) = max\ (Q(u, A), A \in \mathcal{A}_\mathcal{P}, Q(u, A), \in \mathcal{S}_\mathcal{P},).$$

As a measure of the relative behavior of algorithm $A \in \mathcal{A}$ solution on input $u \in \mathcal{I}_\mathcal{P}$ we will use the ration

$$R_\mathcal{P}(A, u) = A(u) / Q_{max}(u),$$ where $Q_{max}(u)$ is the best known solution for the input $u$ or the best solution in the set of solutions $\mathcal{A}_\mathcal{P}(u)$.

From the point of view of the solution search time we can specify additional measure (used time $\mathcal{T}$) and so experiment is *minimization* problem.

1. a measure defined for solved problem $\mathcal{T}$. In our case this is the time used by algorithm.
2. a set of algorithms $\mathcal{A}_\mathcal{T}$ ;
3. a set of possible inputs (instances) $\mathcal{I}_\mathcal{T}$ ;
4. a map $\mathcal{S}_\mathcal{T}$ which for each $u \in \mathcal{I}_\mathcal{T}$ maps a finite set of algorithms $A \in \mathcal{A}_\mathcal{T}$ used times;

The performance of a set of algorithms $\mathcal{A}$ for input $u \in \mathcal{I}_\mathcal{T}$ we define as

$$\mathcal{A}_\mathcal{T}(u) = min\ (T(u, A), A \in \mathcal{A}_\mathcal{T}, T(u, A), \in \mathcal{S}_\mathcal{T},).$$

As a measure of the relative behavior of algorithm $A \in \mathcal{A}$ used time on input $u \in \mathcal{I}_\mathcal{T}$ we will use the ration

$$R_\mathcal{T}(A, u) = T_{min}(u) / \mathcal{A}_\mathcal{T}(u),$$ where $T_{min}(u)$ is the least used time solution for the input $u$.

In our case $\mathcal{A}_\mathcal{T} = \mathcal{A}$ and $\mathcal{I}_\mathcal{P} = \mathcal{I}_\mathcal{T}$.

## *6.2. Algorithms*

All algorithms were implemented by the authors in C++ and compiled with gcc v. 5.3.0 with (-O3) optimization flag. All runs were made on Intel(R) Core(TM) i3-4160 CPU @ 3.60GHz with 2 GB of RAM.

Each algorithm can solve MCP or MIS problem in accordance with solved problem value (MCP or MIS) flag on the input of algorithm.

We do not include reading (generating) the graph and building the adjacency matrix, since these are common to all algorithms. But we include the time to allocate, initialize and destroy the data structures that are specific to each algorithm.

For initiate solution step we implement 4 sequential heuristic algorithms:

1. FV_BIO  –  First Vertex the Best In Old;

The input subgraph is represented as an ordered set of vertices. Scanning through each input subgraph vertices, if vertex adjacent with all solution members then add vertex into solution (E1 algorithm of Jonson[6]).

This kind of algorithms is often called *Greedy* algorithm.

2. SD_WON  –  Smallest Degree the Worst Out New;

Put input subgraph into solution set. While solution is not clique, scanning through each solution vertices, select the min degree vertex, remove it and its neighbourhood(E2 algorithm of Jonson[6]).

This kind of algorithms is often called *the Smallest the Last (SL)* algorithm.

3. SD_ext_WON  –  Smallest Degree(e) Extended the Worst Out New;

Similar to the previous one. But we found that not all maximal clique members occur in the solution, they can be removed on the scanning phase. So we try to extend solution by removed vertices.

4. LD_ BIO  -  Largest Degree the Best In Old;

Order subgraph vertices in no decreasing order (the max degree vertex – the first) and apply to it FV_BIO algorithm.

5. LD_ BIN  -  Largest Degree the Best In New;

While input subgraph is not empty, select max degree vertex, put it into solution and remove it and its neighbourhood from the input subgraph (NMCLIQ algorithm of Tomita[16]).

This kind of algorithms is often called *the Largest the First (LF)* algorithm.

6. LS_1_k<SEQ>  -  (1,$k$)-swap Local Search Algorithm, where <SEQ> means some sequential heuristic(1 - 5 upper algorithms) used on the initialization step. So LS_1_k_LD_BIN means that on the 1st phase we apply LD_BIN algorithm and on the 2nd phase we apply our LS_1_k algorithm to the initial solution.

## 6.3. Instances

### 6.3.1. DIMACS benchmark.

These benchmarks were constructed using the maximum clique instances from the second DIMACS Implementation Challenge [6] which has been used extensively for benchmarking purposes in the recent literature on maximum clique and coloring algorithms. This group contains 199 instances: 66(base group) + 14(additional) + 5(random) graphs + 119(coloring instances).

### 6.3.2. BHOSLIB benchmark

Benchmarks with Hidden Optimum Solutions for Graph Problems (Maximum Clique, Maximum Independent Set, Minimum Vertex Cover and Vertex Coloring) (BHOSLIB)[3]. The maximum independent set benchmark instances are directly transformed from forced feasible SAT benchmarks, with the set of vertices and the set of edges respectively corresponding to the set of variables and the set of binary clauses in SAT instances.

The benchmark clique instances are the complements of above mentioned graph instances. This group contains 41 graphs.

### 6.3.3. Sloane benchmark

This is a collection of graphs arising from coding theory [15].

The benchmark clique instances are the complements of above mentioned graph instances. This group contains 32 graphs.

### 6.3.4. Iovanella benchmark.

They built [8] a collection of instances with $n = 2^i$ for i = 7, 8, ..., 14 i.e., for n ∈ {128, 256, 512, 1024, 2048, 4096, 8192, 16384} and density d = {0.5, 0.9}.

This group contains 16 graphs.

### 6.3.5. Random benchmark

To generate collection of random graphs we use procedure GenRn with input params (n_rpt; n0,nI,nN; d0,dI,dN; seed):

C - is collection of generated Rn,d graphs;

```
while(n_rpt > 0) {
    --n_rpt;
    For(n= n0; n < nN; n+=nI)
        For(d= d0; d < dN; d+=dI)
            C.Add(Generate Rn graph(n, d, seed)); // update seed value;
}
```
Our collections are

  C1 =GenRn (100; 250, 250, 999; 0.1,0.2,0.9);

  C2 = GenRn(50; 1000, 500, 9999; 0.1,0.2,0.9);

  C3 = GenRn (10; 10000, 5000, 50000; 0.1,0.2,0.9);

This a collection of random n = (10000, 15000, …, 45000, 50000) vertex graphs with different density values { 0.1, 0.3, 0.5, 0.7, 0.9 }. We generate this collection especially for algorithms used time estimation. So we have 9 * 5 = 45 graph groups and we generated 5 random graphs for each group.

This group contains ~6500 graphs.

### 6.4. Results

We have a large enough list of instances ~7,000 - graphs has been analyzed! For each graph each algorithm generates as maximal clique and maximal independent set. So we investigated each algorithm on each of ~7000 instances two times(for MCP and for MIS). We joint relative algorithms behavior in the table:

|   | Algorithm name | $R_{\mathcal{P}}(A, u)$ - relative solutions measure | $R_{\mathcal{T}}(A, u)$ relative used time measure |
|---|---|---|---|
| 1 | fvsFV_BIO | 0.770 | 0.706 |
| 2 | fvsLS_1_k_FV_BIO | 0.888 | 0.067 |
| 3 | fvsSD_WON | 0.759 | 0.002 |
| 4 | fvsLS_1_k_SD_WON | 0.935 | 0.002 |
| 5 | fvsSD_WON_f | 0.862 | 0.002 |
| 6 | fvsLS_1_k_SD_WON_f | 0.943 | 0.002 |
| 7 | fvsLD_BIO | 0.830 | 0.303 |

| 8 | fvsLS_1_k_LD_BIO | 0.932 | 0.076 |
| 9 | fvsLD_BIN | 0.980 | 0.113 |
| 10 | fvsLS_1_k_LD_BIN | 0.991 | 0.063 |

The best algorithms from the point of view of algorithms solution deviation from the best solution obtained by algorithms group are fvsLS_1_k_LD_BIN = 0.991 and fvsLD_BIN = 0.980, fvsLS_1_k_SD_WON_f = 0.943, fvsLS_1_k_SD_WON = 0.935.

From the table we see that the fastest algorithm is fvsFV_BIO with relative used time coefficient equals 0.706, then fvsLD_BIO = 0.303, then fvsLD_BIN = 0.113 and so on.

## *Conclusions*

The (1-k)-swap local search fvsLS_1_k_LD_BIN algorithm showed the best solution results but it twice slower than the sequential fvsLD_BIN algorithm .

The sequential fvsLD_BIN algorithm is quite good from the solution point of view and is better than LS algorithms from the used time criteria.

Class of the Smallest Degree the Worst Out algorithms can be removed from the practice on their low speed criteria.